\documentclass[12pt]{amsart}%{article}%
\usepackage{amsmath,amstext,amsfonts,amsthm,amssymb}
\usepackage{eucal}
\usepackage{diagrams}                                      %
\diagramstyle[scriptlabels,PostScript=dvips]               %
\newarrow{Dash}{}{dash}{}{dash}{}                          %
\swapnumbers

\newtheorem{THM}[subsection]{Theorem}

\newtheorem*{Thm}{Theorem}

{  \theoremstyle{definition}

}

\newcommand{\ad}{\mathrm{ad}}

\newcommand{\Bimod}{\operatorname{Bimod}}
\newcommand{\Bicomod}{\operatorname{Bicomod}}

\newcommand{\Lie}{\mathtt{Lie}}

\newcommand{\Coind}{\operatorname{Coind}}
\newcommand{\Com}{\mathtt{Com}}

\newcommand{\Def}{\mathrm{Def}}

\newcommand{\Ext}{\operatorname{Ext}}

\newcommand{\Hom}{\mathrm{Hom}} % straight Hom

\newcommand{\chom}{\mathcal{H}\mathit{om}} % calligraphic Hom for complexes
\newcommand{\Ind}{\operatorname{Ind}}

\newcommand{\Mod}{\mathtt{Mod}}

\newcommand{\Right}{\mathbf{R}}

\newcommand{\Tetra}{\mathtt{Tetra}}

\newcommand{\cE}{\mathcal{E}}

\newcommand{\cP}{\mathcal{P}}

\newcommand{\fg}{\mathfrak{g}}

\newcommand{\F}{\mathbb{F}}

\newcommand{\Q}{\mathbb{Q}}

\newcommand{\wh}{\widehat}

\begin{document}
\title[]{Formality theorem and bialgebra deformations}
\author{Vladimir Hinich}
\address{Department of Mathematics, University of Haifa,
Mount Carmel, Haifa 31905,  Israel}
\email{vhinich@gmail.com}
\author{Dan Lemberg}
\address{Department of Mathematics, University of Haifa,
Mount Carmel, Haifa 31905,  Israel}
\email{lembergdan@gmail.com}
 
\begin{abstract}
In this paper we prove formality of the exterior algebra on $V\oplus V^*$ endowed
with the {\sl big bracket} considered as a graded Poisson algebra. We also discuss
connection of this result to bialgebra deformations of the symmetric algebra of $V$ considered as bialgebra.
\end{abstract}
\maketitle

\section{Introduction}

\subsection{}
\label{ss:H}
In this paper $k$ will always denote a ground field of characteristic zero. 
Fix a finite dimensional vector space $V$ over $k$. This paper deals with the graded
vector space $H=\oplus H^n$ where
\begin{equation}
\label{eq:H}
H^n=\bigoplus_{p+q=n}\wedge^p V\otimes\wedge^q V^*,\quad n\geq 0.
\end{equation}

This vector space has a graded Poisson structure defined as follows. This is a (graded) commutative algebra with unit as  $H=S(W)$ where
$$W=(V\oplus V^*)[-1].$$
Here we use the standard convention for graded versions of commutative algebras, based
on symmetric monoidal structure on the category of graded vector spaces defined by the
commutativity constraint given by the standard formula
\begin{equation}
\label{eq:bracket}
\sigma: X\otimes Y\rTo Y\otimes X,\quad \sigma(x\otimes y)=(-1)^{|x||y|}y\otimes x.
\end{equation}

The commutative algebra $H$ has a degree $-2$ Lie bracket 
(called {\sl big bracket} by Y.~Kosmann-Schwarzbach in~\cite{YKS}) satisfying Leibniz rule with respect to the multiplication, and given on the generators by the formula
\begin{equation}
\label{eq:bracket}
[x,x']=0,\ [y,y']=0,\ [x,y]=[y,x]=\langle y,x\rangle
\end{equation}
for $x,x'\in V$, $y,y'\in V^*$.

Throughout this paper we will be using the language of operads to describe various algebraic structures. Graded Poisson algebras are algebras over a certain operad denoted $P_n$.
Algebras over $P_n$ have a degree zero commutative multiplication and degree $1-n$ Lie bracket
satisfying Leibniz rule. Thus, our algebra $H$ is a $P_3$-algebra.

The operad $P_n$ is Koszul \cite{GJ}, so it has a particularly nice cofibrant resolution 
and a particularly nice notion of homotopy $P_n$-algebra.

In this paper we prove that $H$ is intrinsically formal as $P_3$-algebra. 
This means that any homotopy $P_3$-algebra $X$ such that $H(X)=H$ as $P$-algebras, is 
equivalent to $X$.

The proof follows ideas of Tamarkin \cite{T} and makes use of the criterion
of intrinsic formality described in~\cite{H-Tam}, 4.1.3.

\subsection{} The graded vector space $H$ appears in two interconnected instances in 
deformation theory. The first one is connected to Lie bialgebras, and the second to
associative bialgebras.

\subsubsection{Lie bialgebras}
Recall that, according to Y. Kosmann-Schwarzbach \cite{YKS} a proto-Lie bialgebra structure
on a vector space $V$ is just a degree 3 element $h$ of $H$ satisfying the condition
$[h,h]=0$. Let $\lambda\in V\otimes\wedge^2 V^*$, $\delta\in\wedge^2V\otimes V^*$,
$\alpha\in\wedge^3V$ and $\beta\in\wedge^3V^*$ be the components of $h$. One can easily check 
that if $\alpha=0,\ \beta=0$, the tensors $\lambda$ and $\delta$ determine a Lie bialgebra
structure on $V$. The case $\beta=0$ describes Lie quasi-bialgebras, whereas $\alpha=0$
corresponds to the dual notion.

Lie bialgebras and their "quasi-" versions describe quasiclassical limits of quantized
enveloping algebras as defined by Drinfeld \cite{D}. The notion of Lie proto-bialgebra
naturally generalizes both Lie quasi-bialgebra and coquasi-bialgebras. We are unaware of the existence of the notion of associative proto-bialgebra quantizing Lie proto-bialgebras.

Let $h\in H^3$ satisfy $[h,h]=0$, so that $(V,h)$ is a Lie proto-bialgebra.
The operator $d_h=\ad_h$ is a derivation of both commutative and Lie algebra structure on
$H$, so that $(H,d_h)$ becomes a dg $P_3$-algebra. The dg Lie algebra $(H[2],d_h)$
governs formal deformations of the Lie proto-bialgebra $(V,h)$. In particular, $H[2]$
governs deformations of the commutative Lie bialgebra (in the class of Lie proto-bialgebras).

Our result on formality should be much more relevant to another deformation problem, that of associative bialgebras\footnote{Or associative proto-bialgebras if one could imagine
what they are.}, see \ref{sss:intro-ba}.

\subsubsection{Associative algebras}

Before we start talking about bialgebra deformations, it is worthwhile to remind
what is going on with already classical problem of deformations of associative algebras. 

Let $A$ be an associative algebra. The Hochschild cohomology $HH^*(A)$
has a structure of Gerstenhaber algebra (this is the same as $P_2$-algebra in our notation).
Moreover, $P_2$ is the homology of the small squares operad $E_2$ (F.~Cohen, \cite{C}),
and the Hochschild cochain complex $CC^*(A)$ has a (homotopy)
structure of algebra over $E_2$ (Deligne conjecture). In case $A$ is a polynomial algebra 
its Hochschild cohomology, the algebra of polyvector fields $\wedge T$, is intrinsically formal as $P_2$-algebra, and this result implies
the famous Kontsevich formality theorem for the polynomial ring, see Kontsevich~\cite{KQ}, Tamarkin~\cite{T}. Note that the algebra 
$\wedge T$, similarly to our algebra $H$, can also be interpreted as
the algebra responsible for deformations of the (trivial) Poisson bracket
in the polynomial ring.   

Note that there are two different deformation problems connected to an associative algebra $A$.
Deformations of $A$ are governed by the (shifted and) truncated Hochschild cochain complex
$\Def_A$ concentrated in nonnegative degrees, so that
$$ \Def_A^n=\Hom(A^{\otimes n+1},A),\ n\geq 0.$$
However, this dg Lie algebra is not formal even in case $A$ is a polynomial ring.
The Kontsevich formality theorem states that the full (shifted) Hochschild cochain complex
concentrated in degrees $n\geq -1$ is formal. The deformation problem described by the full Hochschild complex is that of the category of $A$-modules. This deformation problem is not easy to formally define; the difference between two deformation problems can be seen if one studies deformation of a sheaf of associative algebras: as a result of deformation one can get an algebroid stack instead of the deformed sheaf of associative algebras.

\subsubsection{Associative bialgebras} 
\label{sss:intro-ba}
 
Deformation theory for associative bialgebas was pioneered by Gerstenhaber and Schack in \cite{GS} where a deformation complex $C_{GS}(A)$ of a bialgebra $A$ was defined by the ad hoc formulas

\begin{equation}
\label{eq:CGS}
C^n_{GS}(A)=\oplus_{p+q=n}\Hom(A^{\otimes p},A^{\otimes q}),
\end{equation}
with the differential given for $\phi:A^{\otimes p}\to A^{\otimes q}$ by the formula
\begin{equation}
\label{eq:CGSd}
d(\phi)=d_1(\phi)+d_2(\phi),
\end{equation}
where
\begin{eqnarray}
\label{eq:CGSd1}
d_1(\phi)(a_0\otimes\ldots\otimes a_p)=\Delta^{q-1}(a_0)\phi(a_1\otimes\ldots\otimes a_p)+\\
\nonumber+\sum_{i=0}^{p-1}(-1)^{i+1}\phi(a_0\otimes\ldots\otimes(a_ia_{i+1})\otimes\ldots\otimes a_p)+\\
\nonumber+(-1)^{p-1}\phi(a_0\otimes\ldots\otimes a_{p-1})\Delta^{q-1}(a_p)
\end{eqnarray}
and
\begin{eqnarray}
\label{eq:CGSd2}
d_2(a_1\otimes\ldots\otimes a_p)=
(\mu_{p-1}\otimes\phi)\Delta^{\otimes p}(a_1\otimes\ldots\otimes a_p)+\\
\nonumber\sum_{i=1}^q\Delta_i(\phi(a_1\otimes\ldots\otimes a_p))+
(-1)^{p-1}(phi\otimes\mu_{p-1})\Delta^{\otimes p}(a_1\otimes\ldots\otimes a_p).
\end{eqnarray}

Here in formula (\ref{eq:CGSd1}) $\Delta^{q-1}$ denotes the multiple comultiplication $A\rTo A^{\otimes p}$, $\Delta^{\otimes p}$ in formula
(\ref{eq:CGSd2}) denotes the comultiplication induced on $A^{\otimes p}$,
that is, a map $A^{\otimes p}\rTo A^{\otimes p}\otimes A^{\otimes p}$,
$\mu_{p-1}$ denotes the (multiple) product $A^{\otimes p}\rTo A$, and 
$\Delta_i=1\otimes\ldots\otimes\Delta\otimes\ldots\otimes 1$.
Gerstenhaber-Schack cohomology is defined by the formula
$$ H^*_{GS}(A)=H^*(C_{GS}(A)).$$

For the bialgebra $A=S(V)$ with standard multiplication and coproduct, a well-known calculation 
(reproduced below, following Shoikhet~\cite{Sh1}, see~\ref{sss:HSV}) shows $H_{GS}(A)=H$ as a graded vector space. According to a version of Deligne conjecture proven by B.~Shoikhet \cite{Sh2}, the chain complex $C_{GS}$ admits a structure of $E_3$-algebra. Once more, in characteristic zero the operad $E_3$ is equivalent to $P_3$, see~\cite{LV}, so 
Gerstanhaber-Schack cochains admit a canonical homotopy $P_3$-algebra structure.

Therefore, in order to have a complete analog of Kontsevich formality theorem for
bialgebras, one needs to solve two problems.
\begin{itemize}
\item[1.] Verify that the $P_3$-algebra structure on $H$ defined in \ref{ss:H} comes
from the $E_3$-algebra structure on the cochains $C_{GS}(S(V))$.
\item[2.] Understand in what sense $E_3$-algebra $C_{GS}(A)$ governs deformations of 
the bialgebra $A$.
\end{itemize}

We will now describe what we can say about the above problems.

\subsubsection{}
First of all, Gerstenhaber-Schack cohomology is known to be described as $\Ext(A,A)$
calculated in the abelian category of $A$-tetramodules, see Taillefer~\cite{Tf} and 
Section~\ref{sec:HSV} below.
The commutative algebra structure on $H_{GS}(A)$ comes from Yoneda product in Ext's,
and it is not difficult to verify that for $A=S(V)$ this yields the commutative multiplication in $H$ coming from the symmetric algebra structure.

Furthermore, Leibniz rule together with degree considerations imply that the Lie bracket
on the cohomology induced from the $E_3$ structure on the cochains, is determined by its restriction to $W$, that is by a symmetric bilinear form on $V\oplus V^*$. 

One can easily deduce from this that the bracket on $H$ defined by $E_3$ structure is proportional to the one given by the formula (\ref{eq:bracket}). Unfortunately, this does not imply that the bracket is nonzero. This remains a  problem. We were unable to make an explicit computation of the bracket using Shoikhet's description of $E_3$ structure on the chain complex. We believe that a correct way of doing so 
would be using deformation theory; at the moment we are only able to deduce this fact from 
Conjecture 1, see Section~\ref{sec:spec}. 

\subsubsection{}
Gerstenhaber and Shack
used their cohomology to describe obstruction theory: third cohomology of a certain subcomplex
of $C_{GS}(A)$ describes infinitesimal deformations, with obstructions dwelling in the fourth cohomology.  They conjectured the existence of Lie algebra bracket on cohomology so that
the obstruction of infinitesimal deformation given by $u\in H^3_{GS}(A)$ is $[u,u]\in H^4_{GS}(A)$.

Merkulov and Vallette \cite{MV} proved existence of such bracket on
a certain subcomplex of $C_{GS}(A)$; unfortunately, we see no way of
comparing this bracket with the one coming from $E_3$-structure. 

One should also have in mind that, similarly to the case of associative algebras, one cannot expect the  full Gerstenhaber complex to govern deformations of bialgebras; the full complex
should rather govern deformations of a certain "linear" object attached to a bialgebra $A$.

A sensible candidate would be the two-category of categories, left-tensored over the monoidal category of left $A$-modules.

We hope to be able to make sense of this claim in a later publication.

\section{Intrinsic formality of $H$}
\label{sec:intr-form}

The rational homology $P_n$ of the topological operad $E_n$
was calculated by Fred Cohen in 1973. This is a graded
operad over $\Q$ generated by two operations: commutative 
associative multiplication $\mu$ in degree zero, and a Lie bracket $\lambda$ in degree $1-n$, subject to the graded
version of Leibniz rule.

Furthermore, the operad $C_\bullet(E_n,\Q)$ or rational chains is known to be formal: 
it is quasiisomorphic to $P_n$ as an operad of complexes, see \cite{LV}. 

In this section we prove the intrinsic formality of the
$P_3$-algebra $H=S(V[-1]\oplus V^*[-1])$ in the sense of
\cite{H-Tam}, 4.1.2.

We will follow the Tamarkin's idea \cite{T,H-Tam}.

Recall that $\cP=P_3$ is Koszul operad and the free $\cP$-algebra spanned by a 
complex $X$ has the following form.
\begin{equation}
\F_\cP(X)=\F_\Com\circ\F_{\Lie\{2\}}(X)
\end{equation}
where $\F_\Com$ is the free commutative (=symmetric) algebra, $\Lie\{2\}$
is the operad defined by the property that $\Lie\{2\}$-algebra structure on $X$
is the same as $\Lie$-algebra structure on $X[2]$.  

Also $\cP^\perp=\cP\{-3\}^*$, so that the cofree $\cP^\perp$-coalgebra spanned by $H$
has form
\begin{equation}\label{coH}
\F^*_{\cP^\perp}(H)=\F^*_\Com(\F^*_\Lie(H[1])[2])[-3].
\end{equation}
 
We will use the following criterion of intrinsic formality. 
 
\begin{THM}
\label{thm:TAM-criterion}
(see \cite{H-Tam}, 4.1.3).
Let $\fg$ be the dg Lie algebra of coderivations of
$(\F^*_{\cP^\perp}(H),Q)$, where the differential $Q$
is defined by the $\cE$-algebra structure on $H$.
Denote
\begin{equation}
\fg_{\geq 1}=\Hom(\oplus_{i\geq 2}\F^{*i}_{\cP^\perp}(H),H)
\subset \Hom(\oplus_{i\geq 1}\F^{*i}_{\cP^\perp}(H),H)=\fg.
\end{equation}
Then , if the map $H^1(\fg_{\geq 1})\rTo H^1(\fg)$ is zero,
the $\cP$-algebra $H$ is intrinsically  formal. 
\end{THM}

The dg Lie algebra $\fg$ is obtained from a bicomplex,
\begin{equation}
 \fg=\oplus_{p,q\geq 0}\fg^{p,q},\quad
\fg_{\geq 1}=\oplus_{(p,q)\ne(0,0)}\fg^{p,q},
\end{equation}
where 
\begin{equation}
\fg^{p\bullet}=\Hom(\F^{*(p+1)}_\Com(\F^*_\Lie(H[1])[2]),H[3]),
\end{equation}
 and $q+1$ is the total $\Lie$-degree. The horizontal and the vertical components  $Q_l$ and $Q_m$ of the differential 
are defined by the Lie bracket and commutative multiplication on $H$ respectively. In order to calculate the cohomology of $\fg$ we can use the spectral sequence
of the above bicomplex. Look at the complexes
\begin{equation}
(\fg^{p\bullet},Q_m)=\Hom_H(S^{p+1}_H(\F^*_\Lie(H[1])\otimes H[2]),H[3]).
\end{equation} 
The complex $Z:=\F^*_\Lie(H[1])\otimes H$ is the homological Harrison complex of the commutative algebra $H$. Therefore,
\begin{equation}
(\fg^{p\bullet},Q_m)=\Hom_H(S^{p+1}_H(Z[2]),H[3]).
\end{equation}
Now recall that $H=S(W)$ where $W=V[-1]\oplus  V^*[-1]$, so
$Z$ is quasiisomorphic to the shifted module of differentials, $Z=H\otimes W[1]$.

Therefore, the homology of the complex $\fg^{p\bullet}$
with respect to the vertical differential $Q_m$ is
\begin{equation}
E_1^{p,q}=
\begin{cases}
\Hom(S^{p+1}(W[3]),S(W)[3])&\mathrm{ if\ } q=0\\
0 &  \mathrm{otherwise}. 
\end{cases}
\end{equation}

The spaces $E_1^{p,0}$ are quotients of 
$\fg^{p,0}=\Hom(\F^{*p+1}_\Com(H[3]),H[3])$,
so the differential $Q_l$ on $E_0^{p,0}$ is induced by the
Chevalley-Eilenberg differential 
\begin{equation}\label{eq:CE}
(df)(a_1\cdots a_{p+1})=
-\sum_{i=1}^{p+1}(-1)^{i+|a_i|(|a_1|+\ldots+|a_{i-1}|)}[a_i,f(a_1\ldots \wh{a_i}\ldots a_{p+1})]
\end{equation}

The spectral sequence degenerates at term 2 ($E^{p,q}_2=E^{p,q}_\infty$), 
$E_2^{p,q}=0$ for $q>0$ and $E_2^{p,0}$ is the $p$-th cohomology of the complex
$(E_1^{p,0},Q_l)$.

We can now verify the condition of Theorem~\ref{thm:TAM-criterion}.
One-cochain in $E_1^{p,0}$ has form $\sum_p f_p$ with 
$$f_p\in \Hom(S^{p+1}(W[3]),S^{4-2(p+1)}(W)).$$
This immediatelly implies that $f_p=0$ for $p>1$. Such cochain is in the image of $\fg_{\geq 1}$
iff $f_0=0$. Thus, any one-cochain coming from $\fg_{\geq 1}$ is presented by a map
$f_1:S^2(W)\rTo k$.  We will show it is always a boundary. More precisely, we claim
there exists $g\in E^{0,0}_1$ of total degree zero such that $f_1=Q_l(g)$.
The elements of total degree zero in $E_1^{0,0}$ are maps $g:W\rTo W$. The formula (\ref{eq:CE})
shows that for such $g$ its differential is calculated as
$$ Q_l(g)(a,b)=[a,g(b)]+[b,g(a)].$$
Since the bracket restricted to $W$ is a nondegenerate symmetric bilinear form, 
existence of $g$ is a basic fact of linear algebra.
This proves the theorem.
\qed

\section{Gerstenhaber-Schack cohomology of $U\fg$}
\label{sec:HSV}

According to Taillefer~\cite{Tf}, Gerstenhaber-Schack cohomology of a bialgebra $A$
is just $\Ext(A,A)$  in a certain abelian category $\Tetra_A$, the category of $A$-tetramodules. According to Shoikhet~\cite{Sh2}, the braided monoidal structure on the category of tetramodules induces $E_3$-algebra structure on GS cochains. This induces a $P_3$-algebra structure on the cohomology. The graded space $H$ studied in the previous section is the GS cohomology of the
bialgebra $A=S(V)$. We would like to identify the $P_3$-structure on $H$ defined by the formulas (\ref{eq:H})--(\ref{eq:bracket}) with the one induced on $H$ as the GS cohomology. We have not completely succeeded in this.

We check that the commutative multiplication in $H$ induced from the $P_3$ structure comes from the presentation $H=S(W)$; furthermore, we prove that the bracket is proportional
to the one given by formula~(\ref{eq:bracket}). Thus, if the bracket on $H$ is nonzero, the
formality calculation of Section~\ref{sec:intr-form} is applicable. Unfortunately, we were 
unable to prove nonvaishing of the bracket on $H$ induced from the $E_3$-structure on GS cochains. 

We believe that the Lie bracket on $H$ is in fact given by the formula (\ref{eq:bracket}).
We support this belief in Section~\ref{sec:spec} with some speculations and conjectures.
 
\

In this section we present the calculation of $\Right\Hom_{\Tetra_A}(A,A)$
for $A=U\fg$ the enveloping algebra. We present it by a dg algebra which induces Yoneda
product on the cohomology. 

In the special case $\fg=V$ is a commutative Lie algebra,
this allows one to identify the Yoneda product on $H$ with the multiplication in the symmetric algebra. Then we deduce that the bracket is proportional to (\ref{eq:bracket}).

\

The category $\Tetra_A$ of $A$-tetramodules has enough injectives; but it is more convenient 
to make calculations using formalism of $(P,Q)$ pairs described in \cite{Sh1}. In the following subsection we recall the relevant definitions. 

\subsection{Tetramodules}

Let $A$ be a bialgebra.

Recall that a tetramodule structure on a vector space $M$ is a bialgebra structure
on the direct sum $A\oplus M$ sich that the natural projection
$p:A\oplus M\rTo A$ is a bialgebra morphism and $p$ is an abelian group object in the category  of bialgebra morphisms $B\rTo A$ with target $A$. The latter means that the
maps $M\otimes M\to M$ and $M\to M\otimes M$ defined by restriction of the multiplication and the comultiplication on $M$, vanish.

Thus, a tetramodule has both a bimodule and a bicomodule structure, satisfying certain compatibilities. 

The category of $A$-tetramodules is denoted $\Tetra_A$. This is an abelian category with 
enough injectives, see \cite{Tf}. In case $A$ is a Hopf algebra, it is equivalent to the category of Yetter-Drinfeld
modules and is Drinfeld double of the monoidal category of left $A$-modules, see \cite{Schau}.

\subsubsection{Induced and coinduced tetramodules}

 We have two pairs of adjoint functors 

\begin{equation}
\Bicomod_A\pile{\rTo^\Ind\\\lTo_G} \Tetra_A \pile{\rTo^G\\ \lTo_\Coind} \Bimod_A
\end{equation}
where $G$ denotes the forgetful functors, $\Ind$ is the induction and $\Coind$ the coinduction functor defined as in Shoikhet~\cite{Sh1}.

Any tetramodule embeds into a coinduced tetramodule and is an image of induced tetramodule.
Therefore, any tetramodule $X$ admits an induced resolution
$$ \rTo P_n\rTo\ldots\rTo P_0\rTo X\rTo 0$$
and a coinduced resolution
$$ 0\rTo X\rTo Q^0\rTo\ldots\rTo Q^n\rTo\ .$$

One has
\begin{Thm}(see~\cite{Sh1})
One can calculate $\Right\Hom_{\Tetra_A}(X,Y)$ using induced resolution for $X$ and coinduced resolutions for $Y$.
\end{Thm} 

\subsection{The case $A=U\fg$}

Let $A=U\fg$ be the enveloping algebra of a finite dimensional Lie algebra $\fg$
considered as a bialgebra.
 
We will calculate $\Right\Hom_{\Tetra_A}(A,A)$
using an induced and a coinduced resolutions for the tetramodule $A$.
We define $P_n=\Ind(\wedge^n\fg)$ and $Q^n=\Coind(\wedge^n\fg)$,
where in the first formula $\wedge^n\fg$ has the trivial bicomodule structure, whereas in the second formula it has the trivial bimodule structure.

We will use the following notation. For a subset $I\subset N=\{1,\ldots,n\}$ and a collection
of elements $x_i\in\fg,\ i\in N$, we denote as $x_I$ the product $x_{i_1}\wedge\ldots\wedge x_{i_{|I|}}$.

The induced tetramodules $P_n$ form a complex $P_\bullet$ with $H^0(P_\bullet)=A$, with the differentials $\partial_n:P_n\rTo P_{n-1}$ defined  
by the formula
\begin{eqnarray}
\label{eq:d-P}
\partial_n(a\otimes x_N\otimes b)=
\sum_{i=1}^n(-1)^{i-1}ax_i\otimes x_{N-\{i\}}\otimes b-
a\otimes x_{N-\{i\}}\otimes x_ib+\\
\nonumber+\sum_{i<j\in N}(-1)^{i+j}a\otimes [x_i,x_j]\wedge x_{N-\{i,j\}}\otimes b.
\end{eqnarray}
The formulas for a differential in the coinduced resolution $Q^\bullet$ of $A$
are similar: \footnote{But simpler as there is no term coming from (co)bracket.}
the differential $d_n:Q^n\rTo Q^{n+1}$  is given by the formula
\begin{equation}
\label{eq:d_Q}
d_n(a\otimes x_N\otimes b)=\Delta^1_r(a)\otimes\Delta^2_r(a)\wedge x_N\otimes b
- a\otimes x_N\wedge\Delta^1_l(b)\otimes\Delta^2_l(b),
\end{equation}
 where $\Delta^1_r(a)\otimes\Delta^2_r(a)$ (resp., $\Delta^1_l(a)\otimes\Delta^2_l(a)$ )
denotes the projection of $\Delta(a)$ to $U\fg\otimes\fg$ (resp., to $\fg\otimes U\fg$).

\

Now Gerstenhaber-Schack cohomology of $A=U\fg$ can be expressed as
\begin{equation}
\label{eq:GS-pol}
H_{GS}(A)=\chom_{\Tetra_A}(P_\bullet,Q^\bullet),
\end{equation}
where $\chom$ denotes the complex of morphisms in $\Tetra_A$

The right-hand side of the equation can be easily calculated. This is
$H=\wedge\fg\otimes\wedge\fg^*$ as a graded commutative algebra, with the differential 
$d=\ad_\lambda$ where $\lambda\in H^3$ is the tensor defining Lie bracket on $\fg$, and $\ad_\lambda$ makes use of the Lie bracket 
defined by the pairing $\fg^*\otimes \fg\rTo k$.
 
In particular, if $\fg=V$ is a commutative Lie algebra, one has
\begin{equation}
\label{eq:GS-pol}
H_{GS}(A)=\Hom(\wedge V,\wedge V)=H.
\end{equation}
 
We will now describe the Yoneda multiplication in $H_{GS}(U\fg)$.
 
Since tetramodules form a monoidal (even braided monoidal) category, Yoneda product can be expressed via the monoidal stucture as follows.

If $\alpha$ and $\beta$ are cycles in 
$\chom(P_\bullet, Q^\bullet)$
of degrees $m$ and $n$, one has a cycle
$$ \alpha\otimes\beta \in \chom(P_\bullet\otimes P_\bullet,Q^\bullet\otimes Q^\bullet)$$
of degree $m+n$ which yields an element in $H^{m+n}_{GS}(A)$
as $P_\bullet\otimes P_\bullet$ and $Q^\bullet\otimes Q^\bullet$ are also resolutions of $A$.
One can further simplify the formulas using the coalgebra structure on $P_\bullet$ and
the algebra structure on $Q^\bullet$ described as follows.

The forgetful functors $\Tetra_A\rTo\Bicomod_A$ and
$\Tetra_A\rTo\Bimod_A$ are monoidal. Thus, $\Ind$ is colax monoidal functor, that is one has a natural morphism
$$ \Ind(X\otimes Y)\rTo \Ind(X)\otimes\Ind(Y).$$

Similarly, $\Coind$ is lax monoidal, that is one has a canonical morphism
$$ \Coind(X)\otimes\Coind(Y)\rTo\Coind(X\otimes Y).$$

Taking this into account, we can define a quasiisomorphism
$P_\bullet\rTo P_\bullet\otimes P_\bullet$ of complexes of tetramodules as follows.
For $n=p+q$ one has a map $\wedge^nV\to\wedge^pV\otimes\wedge^qV$
(of trivial bicomodules over $SV$) which add up to the commutative comultiplication in the 
algebra $\wedge V$.
This yields the map
$$ P_n\rTo \Ind(\wedge^pV\otimes\wedge^qV)\rTo P_p\otimes P_q.$$

The comultiplication on $P_\bullet$ defined by these maps commutes with the differentials $P_n\rTo P_{n-1}$ defined by the formula (\ref{eq:d-P}). Dually, one has a multiplication 
$Q^\bullet\otimes Q^\bullet\rTo Q^\bullet$ (co)induced by the multiplication in $\wedge V$.

The complex $\chom_{\Tetra_A}(P_\bullet,Q^\bullet)$ has, therefore, a dg commutative algebra
structure which induces the Yoneda product in cohomology.

\subsubsection{$\fg=V$ is commutative}
\label{sss:HSV}

In this case the complex $\chom_{\Tetra_A}(P_\bullet,Q^\bullet)=\Hom(\wedge V,\wedge V)$ has zero differential and is isomorphic to $H$. An easy calculation shows that the 
Yoneda product in this case is simply given by the commutative product in the presentation 
$$H=\wedge V^*\otimes \wedge V.$$
 
\subsection{The Lie bracket}

Here we assume $A=S(V)$.
Since Lie bracket on $H$ should satisfy Leibniz rule, it is uniquely defined by its value on
algebra generators, that is on $W=(V\oplus V^*)[-1]$. Since the bracket has to have degree $-2$
and $H_0=k$, it has to be given by a symmetric bilinear form on $V\oplus V^*$.

Let us show that the bracket has to be proportional to the one defined by the formulas
(\ref{eq:bracket}). In fact, the group $GL(V)$ acts by automorphisms on the bialgebra $A=S(V)$.
Any automorphism $g\in GL(V)$ gives rise to a braided autoequivalence of the category $\Tetra_A$.
Therefore, the (homotopy) $E_3$-algebra structure on $H$ has to be $GL(V)$-equivariant.
But the formula (\ref{eq:bracket}) is the only, up to scalar, $GL(V)$-invariant symmetric
bilinear form on $V\oplus V^*$ as
$$ S^2(V\oplus V^*)=S^2(V)\oplus S^2(V^*)\oplus V\otimes V^*,$$
$S^2(V)$ and $S^2(V^*)$ have no invariants and $V\otimes V^*$ has one-dimensional invariant subspace.

This proves our claim.

\section{Speculations}
\label{sec:spec}

The calculation of Gerstenhaber-Schack cohomology of $U\fg$ presented in the previous section
yields, in particular, a canonical map 
\begin{equation}
\label{eq:lie-to-ba}
\Right\Hom_{U\fg}(k,\fg)\rTo \Right\Hom_{\Tetra_{U\fg}}(U\fg,U\fg)[1]
\end{equation}
(compare to \cite{LM}, Theorem 2$^\prime$).

We want to look at this map as a categorification of an embedding
\begin{equation}
T\rTo \wedge T
\end{equation}
from the Lie algebra of vector fields on a smooth affine variety to the algebra of 
polyvector fields endowed with the Schouten bracket.

Our reasoning is as follows. The left-hand side of the formula, cut and shifted by one,
is the dg Lie algebra governing deformations of Lie algebra $\fg$. The whole shifted left-hand side $\Right\Hom_{U\fg}(k,\fg)[1]$ has a dg Lie algebra structure since it identifies with
the dg Lie algebra of coderivations of the standard Chevalley-Eilenberg chain complex of $\fg$. This is the dg Lie algebra governing deformations of the category of $\fg$-modules {\sl considered as symmetric monoidal category}.

The shifted right-hand side of (\ref{eq:lie-to-ba}), 
$\Right\Hom_{\Tetra_{U\fg}}(U\fg,U\fg)[2]$, is expected to govern deformations of the same
category of $\fg$-modules considered as monoidal category. \footnote{This is not precise.
It should be rather responsible for deformations of 2-category of categories, left-tensored
over $\Mod_\fg$.}

Thus, the map~(\ref{eq:lie-to-ba}) corresponds to the embedding of symmetric monoidal deformations of if the monoidal category $\Mod_\fg$ into its monoidal deformations.

There is no doubt the following claim should be true.

{\bf Conjecture 1.} {\sl The map (\ref{eq:lie-to-ba}) preserves Lie bracket in cohomology,
where Lie bracket in the left-hand side comes from its interpretation as the complex of coderivations, whereas the Lie bracket in the right-hand side is induced from the $E_3$-algebra
structure on the GS cochains.}

This conjecture immediately implies that the Lie bracket in $H=H_{GS}(S(V))$ is in fact given by the formula~(\ref{eq:bracket}). 

Actually, we believe a much stronger conjecture is true.

Recall the version of Kontsevich formality for smooth commutative dg algebras proven in 
\cite{HL}:

\begin{Thm} Let $A$ be a smooth commutative dg algebra over a field of characteristic zero.
Then the Hochschild cochain complex of $A$ is equivalent to the dg algebra of polyvector fields
as (homotopy) Gerstenhaber algebras.
\end{Thm}

We believe that a categorified version of the above result should be valid.

The complex $L=\Right\Hom_{U\fg}(k,\fg)[1]$ is a Lie algebroid over the commutative dg algebra $A=C^*(\fg,k)$. This implies that, similarly to the algebra of polyvector fields
acquiring a Gerstenhaber algebra structure, the shifted symmetric algebra
$S_A(L[-2])=S_A(\Right\Hom_{U\fg}(k,\fg)[-1])$ acquires a structure of $P_3$-algebra. 

Recall that the operads $P_3$ and $E_3$ are equivalent in characteristic zero.

Keeping in mind this equivalence, we believe the following to be true.

{\bf Conjecture 2.}  {\sl The $E_3$-algebra $\Right\Hom_{\Tetra_{U\fg}}(U\fg,U\fg)$ is equivalent to the algebra $S_A(L[-2])$ where $L$ is the dg Lie algebroid
$\Right\Hom_{U\fg}(k,\fg)[1]$ over $C^*(\fg,k)$.
}

\

Our formality result of Section~\ref{sec:intr-form} shows that in case $\fg$ is commutative, Conjecture~1 implies Conjecture~2.


\begin{thebibliography}{MMMM}
\bibitem[C]{C} F.~Cohen, Cohomology of braid spaces, Bulletin AMS, 79 (1973), 763--766.
\bibitem[Cr]{crans} S.~Crans, Quillen closed model structures for sheaves, JPAA, 101, 1995, 35--57.
\bibitem[D]{D} V.~Drinfeld, Quasi-Hopf algebras, Leningrad Math. Journal, 1 (1990), 1419--1457.
\bibitem[GJ]{GJ} E.~Getzler, J.~Jones, Operads, homotopy algebra and iterated integrals for double loop spaces, arXiv:hep-th/9403055.
\bibitem[GS]{GS} M.~Gerstenhaber, S.~Schack, {\sl Bialgebra cohomology, deformations and quantum groups}, Proceedings of NAS 87(1990), 478--481.
\bibitem[H]{H-Tam} V.~Hinich, {\sl Tamarkin's proof of Kontsevich formality theorem}, Forum Mathematicum, 15 (2003), 591--614.
\bibitem[HL]{HL} V.~Hinich, D.~Lemberg, {\sl Noncommutative unfolding of hypersurface singularity}, accepted to J. Noncommutative geometry.
\bibitem[K]{KQ} M.~Kontsevich, {\sl Deformation quantization of Poisson manifolds}, Lett. Math. Phys. 66(2003), 157--216.
\bibitem[KS]{YKS} Y.~Kosmann-Schwarzbach, {\sl
Grand crochet, crochets de Schouten et cohomologies d'algèbres de Lie},   C. R. Acad. Sci. Paris Sér. I Math. 312 (1991), no. 1, 123–126.
\bibitem[LV]{LV} P.~Lambrechts, I.~Volic, {\sl Formality of the little N-discs operad}, Memoirs AMS, 230 (2014), no. 1079.
\bibitem[L:HA]{L-HA} J.~Lurie, {\sl Higher algebra}, manuscript available from the author's homepage   (http://www.math.harvard.edu/~lurie/).
\bibitem[L:DAGX]{L-DAGX} J.~Lurie, {\sl Formal moduli problems (DAG X)},
manuscript available from the author's homepage (http://www.math.harvard.edu/~lurie/). 
\bibitem[L:M]{L.M} J. Lurie, {\sl Moduli problems for ring spectra},
manuscript available from the author's homepage (http://www.math.harvard.edu/~lurie/).
\bibitem[LM]{LM} A.~Lazarev, M.~Movshev, {\sl Deformations of Hopf algebras}, Russian Math. Surveys, 46(1991), 253--254.
\bibitem[MCS]{mccs} J.~McClure, J.~Smith, A solution of Deligne's Hochschild cohomology conjecture
\bibitem[MV]{MV} S.~Merkulov, B.~Vallette, Deformation theory of representations of prop(erad)s. I. J. Reine Angew. Math., 634(2009), 
51--106.
\bibitem[Schau]{Schau}P.~Schauenburg, {\sl Hopf modules and Yetter-Drinfeld modules}, J. Algebra 169 (1994), 874--890.
 \bibitem[Sh1]{Sh1} B.~Shoikhet, Tetramodules over a bialgebra form a
2-fold monoidal category, Appl. Categ. Structures, 21 (2013), 291--309.
\bibitem[Sh2]{Sh2} B.~Shoikhet, Differential graded categories and Deligne conjecture, arXiv 1303.2500.
\bibitem[Tf]{Tf} R.~Taillefer, Injective Hopf bimodules, cohomologies of infinite-dimensional Hopf algebras and graded commutativity of the Yoneda product, J. Algebra 276 (2004) 259--279
\bibitem[T]{T} D.~Tamarkin, Another proof of Kontsevich formality theorem, arXiv:math/9802...
 
\end{thebibliography}
\end{document}